\documentclass[10pt,a4paper]{article}

\usepackage{amssymb}
\usepackage{amsmath}
\usepackage{amsthm}
\usepackage[top=1.5in,bottom=1.5in]{geometry}

\newtheorem{theorem}{Theorem}[section]
\newtheorem{proposition}[theorem]{Proposition}
\newtheorem{corollary}[theorem]{Corollary}

\newtheorem{remark}[theorem]{Remark}
\newtheorem{lemma}[theorem]{Lemma}

\begin{document}
\title{Mutually nearest and mutually farthest points of sets in geodesic spaces}
\author{Rafa Esp\' inola$^{*}$, Adriana Nicolae$^{**}$}
\date{}
\maketitle

\begin{center}
{\footnotesize
$^{*}$Dpto. de An\'alisis Matem\'atico\\Universidad de Sevilla\\Apdo. 1160, 41080 Sevilla, Spain\\(email: espinola@us.es)\\
\ \\
$^{**}$Department of Applied Mathematics\\Babe\c s-Bolyai University\\ Kog\u alniceanu 1, 400084, Cluj-Napoca, Romania\\ (email: anicolae@math.ubbcluj.ro)
}
\end{center}

\begin{abstract}
Let $A$ and $X$ be nonempty, bounded and closed subsets of a geodesic metric space $(E,d)$. The minimization (resp. maximization) problem denoted by $\min(A,X)$ (resp. $\max(A,X)$) consists in finding $(a_0,x_0) \in A \times X$ such that $d(a_0,x_0) = \inf\left\{d(a,x) : a \in A, x \in X\right\}$ (resp. $d(a_0,x_0) = \sup\left\{d(a,x) : a \in A, x \in X\right\}$). We study the well-posedness of these problems in different geodesic spaces considering the set $A$ fixed. Let $P_{b,cl,cv}(E)$ be the space of all nonempty, bounded, closed and convex subsets of $E$ endowed with the Pompeiu-Hausdorff distance. We show that in a space with a convex metric, curvature bounded below and the geodesic extension property, the family of sets in $P_{b,cl,cv}(E)$ for which $\max(A,X)$ is well-posed is a dense $G_\delta$-set in $P_{b,cl,cv}(E)$. We give a similar result for $\min(A,X)$ without needing the geodesic extension property. Besides, we analyze the situations when one set or both sets are compact and prove some results specific to CAT$(0)$ spaces. We also prove a variant of the Drop theorem in geodesic spaces with a convex metric and apply it to obtain an optimization result for convex functions.
\end{abstract}

\section{Introduction}
Let $(E,d)$ be a metric space, $A \subseteq E$ nonempty and closed (resp. nonempty, bounded and closed), and $x \in E \setminus A$. The nearest point problem (resp. farthest point problem) of $x$  to $A$ consists in finding a point $a_0 \in A$ (the solution of the problem) such that $d(x,a_0) = \inf\{d(x,a) : a \in A\}$ (resp. $d(x,a_0)= \sup\{d(x,a) : a \in A\}$). Ste\v ckin \cite{Steckin} was one of the first who realized that in case $E$ is a Banach space, the geometric properties like strict convexity, uniform convexity, reflexivity and others play an important role in the study of nearest and farthest point problems. His work triggered a series of results so-called ``in the spirit of Ste\v ckin" because the ideas he used were adapted again and again by different authors to various contexts (see \cite{DeBlasi2, DeBlasi, Li}). In \cite{Steckin}, Ste\v ckin proved, in particular, that for each nonempty and closed subset $A$ of a uniformly convex Banach space, the complement of the set of all points $x \in E$ for which the nearest point problem of $x$ to $A$ has a unique solution is of first Baire category. One of the results also given in \cite{Steckin} and later improved by De Blasi, Myjak and Papini in \cite{DeBlasi2} was going to become a key tool in proving best approximation results and was called Ste\v ckin's Lemma. 

In \cite{DeBlasi}, De Blasi, Myjak and Papini studied more general problems than the ones of nearest and farthest points. Namely, they considered the problem of finding two points which minimize (resp. maximize) the distance between two subsets of a Banach space. They focused on the well-posedness of the problem which consists in showing the uniqueness of the solution and that any approximating sequence of the problem must actually converge to the solution (see section 2 for details). The authors proved that if $A$ is a nonempty, bounded and closed subset of a uniformly convex Banach space $E$, the family of sets $X \in P_{b,cl,cv}(E)$ for which the maximization problem $\max(A,X)$ is well-posed is a dense $G_\delta$-set in $P_{b,cl,cv}(E)$, where $P_{b,cl,cv}(E)$ is endowed with the Pompeiu-Hausdorff distance. For the minimization problem $\min(A,X)$ a similar result is proved where $X$ belongs to a particular subspace of $P_{b,cl,cv}(E)$. A nice synthesis of issues concerning nearest and farthest point problems in connection with the geometric properties of Banach spaces and some extensions of these problems can be found in \cite{Cobzas}.

Zamfirescu initiated in \cite{Zam} the investigation of this kind of problems in the context of geodesic spaces. Later on, researchers have focused on adapting the ideas of Ste\v ckin \cite{Steckin} into the geodesic setting. In particular, Zamfirescu proved in \cite{Zamfirescu} that in a geodesic space $E$ without bifurcating geodesics, for a fixed compact set $A$, the set of points $x \in E$ for which the nearest point problem of $x$ to $A$ has a single solution is a set of second Baire category. Motivated by this result, Kaewcharoen and Kirk \cite{Kae&Kirk} showed that if $E$ is a CAT$(0)$ space with the geodesic extension property and with curvature bounded below, for any fixed closed set $A$, the set of points $x \in E$ for which the nearest point problem of $x$ to $A$ has a unique solution is a set of second Baire category. A similar result is proved for the farthest point problem. Very recent results in the context of a space with curvature bounded below were obtained in \cite{Espinola} where the authors prove a variant of Ste\v ckin's Lemma that allows them to give some porosity theorems which are stronger results than the ones in \cite{Kae&Kirk}. 

In this paper we will also be concerned with the geometric result known as the Drop Theorem. The original version of this theorem was proved by Dane\v s \cite{Danes} and is a very useful tool in nonlinear analysis because of its equivalence to the Ekeland Variational Principle. Penot \cite{Penot} proved that in fact, it is also equivalent to the Flower Petal Theorem. In \cite{Georgiev}, generalized versions of the Drop Theorem are proved and afterwards used in the proofs of various minimization problems. For more details see also \cite{Hyers}.

The purpose of this paper is to study in the context of geodesic metric spaces the problem of minimizing (resp. maximizing) the distance between two sets, originally considered by De Blasi, Myjak and Papini in \cite{DeBlasi} for uniformly convex Banach spaces. The given
results rely on a property of the convex hull of a convex set with a point in spaces with convex metric, Lemma \ref{Lemma_ineq_co}, which is proved at the beginning of section \ref{Results_curv_b_b}. We show that if $E$ is a geodesic space with convex metric, curvature bounded below and the geodesic extension property, the family of sets in $P_{b,cl,cv}(E)$ for which $\max(A,X)$ is well-posed is a dense $G_\delta$-set in $P_{b,cl,cv}(E)$. A similar result is given for the minimizing problem, $\min(A,X)$, with no need of the geodesic extension property. These results give natural counterparts to those obtained by De Blasi {\it et al.} in \cite{DeBlasi}. After this we focus on the case of CAT$(0)$ spaces, where the rich geometry of these spaces will be used to relax certain conditions in relation to the well-posedness problem. Then, in section 4, we show that the boundedness condition on the curvature of the space is no longer needed if we impose compactness conditions on the sets. Both minimization and maximization problems are discussed in this context where we replace the condition on the curvature by that of not having bifurcating geodesics introduced by Zamfirescu in \cite{Zamfirescu}. Finally, in our last section, we consider the Drop Theorem in geodesic spaces. With the aid of the Strong Flower Petal Theorem we derive a version of the Drop Theorem in our context which will be used to obtain an optimization result for convex and continuous real-valued functions defined on geodesic spaces.

\section{Preliminaries} \label{Preliminaries}
Let $(E,d)$ be a metric space. A {\it geodesic} in $E$ is an isometry from $\mathbb{R}$ into $E$ (we may also refer to the image of this isometry as a geodesic). A {\it geodesic path} from $x$ to $y$ is a mapping $c:[0,l] \to E$, where $[0,l] \subseteq \mathbb{R}$, such that $c(0) = x, c(l) = y$ and $d\left(c(t),c(t')\right) = \left|t - t'\right|$ for every $t,t' \in [0,l]$. The image $c\left([0,l]\right)$ of $c$ forms a {\it geodesic segment} which joins $x$ and $y$ and is not necessarily unique. If no confusion arises, we will use  $[x,y]$ to denote a geodesic segment joining $x$ and $y$.  $(E,d)$ is a {\it geodesic space} if every two points $x,y \in E$ can be joined by a geodesic path. A point $z\in E$ belongs to the geodesic segment $[x,y]$ if and only if there exists $t\in [0,1]$ such that $d(z,x)= td(x,y)$ and $d(z,y)=(1-t)d(x,y)$, and we will write $z=(1-t)x+ty$ for simplicity. $(E,d)$ has the {\it geodesic extension property} if each geodesic segment is contained in a geodesic. For a very comprehensive treatment of geodesic metric spaces the reader may check \cite{Bridson}.

In a geodesic space $(E,d)$, a {\it function $f : E \to \mathbb{R}$ is convex} if for any geodesic path $c : [0,l] \to E$ we have
\[f(c(tl)) \le (1-t)f(c(0)) + tf(c(l)) \mbox{ for all } t \in [0,1].\]
The {\it metric $d : E \times E \to \mathbb{R}$ is convex} if given any pair of geodesic paths $c_1 : [0, l_1] \to E$ and $c_2 : [0,l_2] \to E$ with $c_1(0) = c_2(0)$ one has 
\[d(c_1(tl_1),c_2(tl_2)) \le td(c_1(l_1),c_2(l_2)) \mbox{ for all } t \in [0,1].\] 
Applying a simple reasoning we notice that we can renounce to the condition $c_1(0) = c_2(0)$. Then,
\[d(c_1(tl_1),c_2(tl_2)) \le (1-t)d(c_1(0),c_2(0)) + td(c_1(l_1),c_2(l_2)) \mbox{ for all } t \in [0,1].\] 
A geodesic space having the metric convex will be referred as a space with convex metric.

A {\it subset $X$ of $E$ is convex} if any geodesic segment that joins every two points of $X$ is contained in $X$. Let $G_1(X)$ denote the union of all geodesics segments with endpoints in $X$. Notice that $X$ is convex if and only if $G_1(X) = X$. Recursively, for $n \ge 2$ we set $G_n(X) = G_1(G_{n-1}(X))$. Then the {\it convex hull} of $X$ will be
\[\mbox{co}(X) = \bigcup_{n \in \mathbb{N}}G_n(X).\]
By $\overline{\mbox{co}}(X)$ we shall denote the closure of the convex hull. It is easy to see that in a geodesic space with convex metric, the closure of the convex hull will be convex and hence it is the smallest closed convex set containing $X$.

Let $\kappa \in \mathbb{R}$ and $n \in \mathbb{N}$. The classical model spaces $M^n_{\kappa}$ are defined in the following way: if $\kappa > 0$, $M^n_{\kappa}$ is obtained from the spherical space $\mathbb{S}^n$ by multiplying the spherical distance with $1/\sqrt{\kappa}$; if $\kappa = 0$, $M^n_{0}$ is the $n$-dimensional Euclidean space $\mathbb{R}^n$; and if $\kappa < 0$, $M^n_{\kappa}$ is obtained from the hyperbolic space $\mathbb{H}^n$ by multiplying the hyperbolic distance with $1/\sqrt{-\kappa}$. For more details about these spaces one can consult \cite{Bridson}.

A {\it geodesic triangle} $\Delta(x_1,x_2,x_3)$ consists of three points $x_1, x_2$ and $x_3$ in $X$ (the {\it vertices} of the triangle) and three geodesic segments corresponding to each pair of points (the {\it edges} of the triangle). For the geodesic traingle $\Delta$=$\Delta(x_1,x_2,x_3)$, a {\it $\kappa$-comparison triangle} is a triangle $\bar{\Delta} = \Delta(\bar{x}_1, \bar{x}_2, \bar{x}_3)$ in $M^2_{\kappa}$ such that $d(x_i,x_j) = d_{M^2_{\kappa}}(\bar{x}_i,\bar{x}_j)$ for $i,j \in \{1,2,3\}$. For $\kappa$ fixed, $\kappa$-comparison triangles of geodesic triangles (having perimeter less than $2\pi/\sqrt{\kappa}$ if $\kappa > 0$) always exist and are unique up to an isometry (see \cite[Lemma 2.14]{Bridson}). 

A geodesic triangle $\Delta$ satisfies the {\it CAT$(\kappa)$} (resp. {\it reversed CAT$(\kappa)$}) {\it inequality} if for every $\kappa$-comparison triangle $\bar{\Delta}$ of $\Delta$ and for every $x,y \in \Delta$ we have
\[d(x,y) \le d_{M^2_{\kappa}}(\bar{x},\bar{y}) \mbox{ (resp. } d(x,y) \ge d_{M^2_{\kappa}}(\bar{x},\bar{y})\mbox{)},\]
where  $\bar{x},\bar{y} \in \bar{\Delta}$ are the corresponding points of $x$ and $y$, i.e., if $x = (1-t)x_i + tx_j$ then $\bar{x} = (1-t)\bar{x}_i + t\bar{x}_j$.

If $\kappa \le 0$, a {\it CAT$(\kappa)$ space} (also known as a space of bounded curvature in the sense of Gromov) is a geodesic space for which every geodesic triangle satisfies the CAT$(\kappa)$ inequality. 
 
A geodesic metric space is said to have {\it curvature bounded below} if there exists $\kappa < 0$ such that every geodesic triangle satisfies the reversed CAT$(\kappa)$ inequality. Other properties of spaces with curvature bounded below and equivalent definitions can be found in \cite{Burago}.

A metric space is said to be {\it without bifurcating geodesics} (see \cite{Zamfirescu}) if for any two segments with the same initial point and having
another common point (different to the initial one), this second point is a common endpoint of both or one segment contains the other. From the definitions, it is easy to see that a space with curvature bounded below cannot have bifurcating geodesics.

CAT$(0)$ spaces are a particular class of CAT$(\kappa)$ spaces which has called the attention of a large number of researchers in the last decades due to its rich geometry and relevance in different problems. The fact that the metric on a CAT$(0)$ space is convex has a great impact on the geometry of the space, but we must mention that having the metric convex is a weaker property than being CAT$(0)$. One can view CAT$(0)$ spaces as variants of Hilbert spaces in the metric setting. In a CAT$(0)$ space we have the following inequality in relation to the generalized parellelogram law of a Hilbert space. Let $x,y_1,y_2$ be points in a CAT$(0)$ space and let $m = (1-t)y_1 +ty_2$ for some $t \in [0,1]$. Then
\[d\left(x,m\right)^2 \le (1-t)d(x,y_1)^2 + td(x,y_2)^2 - t(1-t)d(y_1,y_2)^2.\] 
If above $t = 1/2$, then the inequality is known as the {\it (CN) inequality of Bruhat and Tits} \cite{Bruhat&Tits}. In fact, this inequality is equivalent to the CAT$(0)$ condition. We will refer to the above inequality as the generalized (CN) inequality.\\
Let $(E,d)$ be a complete CAT$(0)$ space and let $(x_n)_{n \in \mathbb{N}}$ be a bounded sequence in $E$. For $x \in E$, define $r(x,(x_n)) = \limsup_{n \to \infty}d(x,x_n)$. The {\it asymptotic radius} of $(x_n)_{n \in \mathbb{N}}$ is given by
\[r((x_n)) = \inf\left\{r(x,(x_n)) : x \in X\right\},\]
and the {\it asymptotic center} of $(x_n)_{n \in \mathbb{N}}$ is the set
\[A((x_n)) = \left\{x \in X : r(x,(x_n)) = r((x_n))\right\}.\]
In \cite{Dhomp&Kirk} it is shown that in a complete CAT$(0)$ space or in a closed convex subset of a complete CAT$(0)$ space, the asymptotic center of a bounded sequence is a singleton.\\
A sequence $(x_n)_{n \in \mathbb{N}}$ in a complete CAT$(0)$ space $E$ is said to {\it $\Delta$-converge} to $x \in E$ if $x$ is the unique asymptotic center of $(u_n)_{n \in \mathbb{N}}$ for every subsequence $(u_n)_{n \in \mathbb{N}}$ of $(x_n)_{n \in \mathbb{N}}$. In this case $x$ will be called the {\it $\Delta$-limit} of $(x_n)_{n \in \mathbb{N}}$ and we will write $\Delta$-$\lim_{n \to \infty}x_n= x$. If $(x_n)_{n \in \mathbb{N}}$ is a bounded sequence in $E$ that $\Delta$-converges to $x$ and if $y \in X$ with $y \ne x$ then, see \cite{Kirk&Panyanak} for detalis,
\[\limsup_{n \to \infty}d(x_n,x) < \limsup_{n \to \infty}d(x_n,y).\]
Every bounded, closed and convex subset of $E$ contains all the $\Delta$-limits of all its $\Delta$-convergent sequences and every bounded sequence in $E$ contains a $\Delta$-convergent subsequence. Based on the stated properties, it is easy to see that a bounded sequence $\Delta$-converges to $x \in E$ provided all its $\Delta$-convergent subsequences have the same $\Delta$-limit $x$.
For more details about the concept of $\Delta$-convergence in CAT$(0)$ spaces one can consult \cite{Kirk&Panyanak}.

We say that the geodesic space $(E,d)$ is {\it reflexive} if every descending sequence of nonempty, bounded, closed and convex subsets of $E$ has nonempty intersection. A simple example of a reflexive metric space is a reflexive Banach space. Other examples include complete CAT$(0)$ spaces, complete uniformly convex hyperbolic spaces with a monotone or a lower semi-continuous from the right modulus of uniform convexity (see \cite{Bozena, Kohl&Leustean}) and others.

Let $(E,d)$ be a metric space. Taking $z \in E$ and $r > 0$ we denote the open (resp. closed) ball centered at $z$ with radius $r$ by $B(z,r)$ (resp. $\widetilde{B}(z,r)$).
Given $X$ a nonempty subset of $E$, we define the {\it distance of a point} $z \in E$ to $X$ by $\mbox{dist}(z,X) = \inf\{d(z,x) : x \in X\}.$ The {\it metric projection} (or {\it nearest point mapping}) $P_X$  onto $X$ is the mapping
\[P_X(y)=\{ x \in X : d(x,y)=\mbox{dist}(y,X)\}, \mbox{ for every } y\in E.\]
The closure of the set $X$ will be denoted as $\overline{X}$.\\
If $X$ is additionally bounded, the {\it diameter} of $X$ is given by $\mbox{diam}X = \sup\{d(x,y) : x , y \in X\}$ and the {\it remotal distance of a point} $z \in E$ to $X$ is defined by $\mbox{Dist}(z,X) = \sup\{d(z,x) : x \in X\}.$ The {\it farthest point mapping} $F_X$ onto $X$ is given by
\[F_X(y)=\{ x \in X : d(x,y)=\mbox{Dist}(y,X)\},\mbox{ for every } y \in E.\]
If $X$ and $Y$ are nonempty and bounded subsets of $E$, one defines the {\it Pompeiu-Hausdorff distance} as
\[h(X,Y) = \max\left\{\sup\{\mbox{dist}(x,Y) : x \in X\}, \sup\{\mbox{dist}(y,X): y \in Y\}\right\}.\]

From now on, if nothing else is mentioned, $E$ will stand for a geodesic metric space. We consider the following families of sets
\[P_{cl}(E) = \left\{X \subseteq E : X \mbox{ is nonempty and closed}\right\},\]
\[P_{b,cl}(E) = \left\{X \subseteq E : X \mbox{ is nonempty, bounded and closed}\right\},\]
\[P_{b,cl,cv}(E) = \left\{X \subseteq E : X \mbox{ is nonempty, bounded, closed and convex}\right\},\]
\[P_{cp}(E) = \left\{X \subseteq E : X \mbox{ is nonempty and compact}\right\},\]
\[P_{cp,cv}(E) = \left\{X \subseteq E : X \mbox{ is nonempty, compact and convex}\right\}.\]
If $E$ is complete, then $P_{b,cl}(E)$ and $P_{cp}(E)$ are complete under the Pompeiu-Hausdorff distance. If, additionally, the metric of $E$ is convex, then, by an easy adaptation of the argument in the Banach space context, one can prove that $P_{b,cl,cv}(E)$ and $P_{cp,cv}(E)$ are also complete with respect to the Pompeiu-Hausdorff distance.

Following \cite{DeBlasi}, for $X,Y \in P_{b,cl}(E)$ and $\sigma > 0$, we set
\[\lambda_{XY} = \inf\left\{d(x,y) : x \in X, y \in Y\right\} \mbox{, } \mu_{XY} = \sup\left\{d(x,y) : x \in X, y \in Y\right\},\]
\[L_{XY}(\sigma) = \left\{x \in X : \mbox{dist}(x,Y) \le \lambda_{XY} + \sigma\right\},\]
\[M_{XY}(\sigma) = \left\{x \in X : \mbox{Dist}(x,Y) \ge \mu_{XY} - \sigma\right\}.\]
The {\it minimization} (resp. {\it maximization}) {\it problem} denoted by $\min(X,Y)$ (resp. $\max(X,Y)$) consists in finding $(x_0,y_0) \in X \times Y$ (the {\it solution} of the problem) such that $d(x_0,y_0) = \lambda_{XY}$ (resp. $d(x_0,y_0) = \mu_{XY}$). A sequence $(x_n,y_n)_{n \in \mathbb{N}}$ in $X \times Y$  such that $d(x_n,y_n) \to \lambda_{XY}$ (resp. $d(x_n,y_n) \to \mu_{XY}$) is called a {\it minimizing} (resp. {\it maximizing}) {\it sequence}. The problem $\min(X,Y)$ (resp. $\max(X,Y)$) is said to be {\it well-posed} if it has a unique solution $(x_0,y_0) \in X \times Y$ and for every minimizing (resp. maximizing) sequence $(x_n,y_n)_{n \in \mathbb{N}}$ we have $x_n \to x_0$ and $y_n \to y_0$. In the following we give a characterization of the well-posedness of $\min(X,Y)$ (resp. $\max(X,Y)$) which can be proved by a straightforward verification of the above definitions. 
\begin{proposition} \label{Prop_well-posed}
Let $(E,d)$ be a complete geodesic metric space and $X,Y \in P_{b,cl}(E)$. The problem $\min(X,Y)$ (resp. $\max(X,Y)$) is well-posed if and only if
\[\inf_{\sigma > 0}\emph{diam}L_{XY}(\sigma) = 0 \mbox{ and } \inf_{\sigma > 0}\emph{diam}L_{YX}(\sigma) = 0,\]
\[(\mbox{resp.} \inf_{\sigma > 0}\emph{diam}M_{XY}(\sigma) = 0 \mbox{ and } \inf_{\sigma > 0}\emph{diam}M_{YX}(\sigma) = 0).\]
\end{proposition}

We give next some results obtained in \cite{Espinola} that will constitute key tools in proving our results.
Let 
\begin{equation}\label{assumption1}
x\in E, r>0, y\in B(x,r/2)\setminus\{x\} \mbox{ and } 0 \le \sigma \le 2d(x,y).
\end{equation}
Set
\[ D(x,y;r,\sigma)=\widetilde{B}(y, r-d(x,y)+\sigma) \setminus B(x,r).\]
Following \cite{Espinola}, for $\kappa \in (-\infty,0)$, define the real function $F_\kappa$ on $\mathbb{R}^3_+$ by
$$ \begin{array}{ll}
F_\kappa(d,r,\sigma)=&
\frac{2}{\sqrt{-\kappa}}\,{\rm arccosh} \bigg(\cosh^2(\sqrt{-\kappa}\,(r-d+\sigma))
  -\frac{\displaystyle\sinh (\sqrt{-\kappa}\,(r-d+\sigma))}{\displaystyle\sinh (\sqrt{-\kappa}\,d)}\\
  &\cdot\left[\cosh (\sqrt{-\kappa}\,r)
- \cosh (\sqrt{-\kappa}\,d)\cosh (\sqrt{-\kappa}\,(r-d+\sigma))\right]\bigg)
\end{array}$$
for each $(d,r,\sigma)\in \mathbb{R}^3_+$.\\
In \cite{Espinola}, the authors prove the following properties of the function $F_\kappa$ and give an estimation of the diameter of the sets $D(x,y;r,\sigma)$. This estimation yields a variant of Ste\v ckin's Lemma for spaces of curvature bounded below.
\begin{proposition} \label{Prop_F_kappa}
The function $F_\kappa$ is continuous on $\mathbb{R}_+^3$ and for any $d\ge0 $ and $r\ge0$, we have that $F_\kappa(d,r,0)=0$.
\end{proposition}
\begin{proposition} \label{Prop_diam}
Let  $(E,d)$ be a geodesic space of curvature bounded below by $\kappa$ and let $x,y,r$ and $\sigma$ satisfy \eqref{assumption1}. Suppose there exists $u\in E$ in a geodesic passing through $x$ and $y$ such that $d(x,u) = r$ and $d(y,u) = r - d(x,y)$. Then the following estimate holds:
\[\emph{diam}D(x,y;r,\sigma) \le F_\kappa (d(x,y),r,\sigma) + 2\sigma .\]
\end{proposition} 

\section{Results in spaces with convex metric and curvature bounded below} \label{Results_curv_b_b}
We begin this section by giving an estimation for $\mbox{dist}(y,X)$, where $X \in P_{b,cv}(E)$, $x' \in E$ such that $\mbox{dist}(x',X) > 0$ and $y \in \overline{\mbox{co}}\left(X \cup \{x'\}\right)$. It is easy to see that in a geodesic metric space with convex metric, $\mbox{dist}(y,X) < \mbox{dist}(x',X)$ for every $y \in \mbox{co}\left(X \cup \{x'\}\right)$ with $y \ne x'$. We sharpen this upper bound in the following way.

\begin{lemma} \label{Lemma_ineq_co}
Let $E$ be a geodesic metric space with convex metric and let $X \in P_{b,cv}(E)$. Suppose $x' \in E$ such that $\emph{dist}(x',X) >0$. Then, for every $y \in \overline{\emph{co}}\left(X \cup \{x'\}\right)$,
\begin{equation}\label{ineq_co}
\emph{dist}(y,X) \le \emph{dist}(x',X) - \frac{\emph{dist}(x',X)}{\emph{dist}(x',X) + \emph{diam}X}d(x',y).
\end{equation}
\end{lemma}
\begin{proof} It is enough to prove (\ref{ineq_co}) for $y \in \mbox{co}\left(X \cup \{x'\}\right)$. For simplicity, let
\[\alpha = \frac{\mbox{dist}(x',X)}{\mbox{dist}(x',X) + \mbox{diam}X}.\]
Let us first prove the inequality for all points belonging to the geodesic segments in $G_1\left(X \cup \{x'\}\right)$ and then we proceed by induction on $G_n\left(X \cup \{x'\}\right)$. Notice that if a geodesic segment has both endpoints in $X$ then the result trivially holds for any of its points. Suppose $x \in X$ and $z_t = (1-t)x + tx'$ for some $t \in [0,1]$. For $\epsilon > 0$ there exists $x^* \in X$ such that $d(x',x^*) < $ dist$(x',X) + \epsilon$. Let $x_t = (1-t)x + tx^*$. Then,
\begin{align*}
\mbox{dist}(z_t,X) \le d(z_t,x_t) \le td(x',x^*) < t\mbox{dist}(x',X) + t\epsilon = \mbox{dist}(x',X) - (1-t)\mbox{dist}(x',X) + t\epsilon.
\end{align*}
Since $d(x',x) \le \mbox{dist}(x',X) + \mbox{diam}X$,
\begin{align*}
\mbox{dist}(z_t,X) \le \mbox{dist}(x',X) - \alpha(1-t)d(x',x) + t\epsilon = \mbox{dist}(x',X) - \alpha d(x',z_t) + t\epsilon.
\end{align*}
Letting $\epsilon \searrow 0$, we obtain the desired inequality.\\
Suppose $(\ref{ineq_co})$ holds for every $y \in G_n\left(X \cup \{x'\}\right)$. We show that it also holds for every $y \in G_{n+1}\left(X \cup \{x'\}\right)$. Take $z_1, z_2\in G_n\left(X \cup \{x'\}\right)$ and let $z_t = (1-t)z_1 + tz_2.$ For $\epsilon > 0$ there exist $x_1,x_2 \in X$ such that \[d(x_1,z_1) < \mbox{dist}(z_1,X) + \epsilon \mbox{ and } d(x_2,z_2) < \mbox{dist}(z_2,X) + \epsilon.\] 
Let $x_t = (1-t)x_1 + tx_2$. Then,
\begin{align*}
\mbox{dist}(z_t,X) &\le d(z_t,x_t) \le (1-t)d(x_1,z_1) + td(x_2,z_2) \\
& <  (1-t)\mbox{dist}(z_1,X) + t\mbox{dist}(z_2,X) + \epsilon\\
& \le \mbox{dist}(x',X) - \alpha\left((1-t)d(x',z_1) + td(x',z_2)\right) + \epsilon\\
& \le \mbox{dist}(x',X) - \alpha d(x',z_t) + \epsilon.
\end{align*} 
Now we only need to let $\epsilon \searrow 0$. Hence, the induction is complete and the conclusion follows.
\end{proof}

The following lemma is an analogue in the metric setting of a property of Banach spaces \cite[Proposition 2.3]{DeBlasi}. 

\begin{lemma} \label{Lemma_conv_diam}
Let $E$ be a geodesic metric space with convex metric and let $X \in P_{b,cv}(E)$. For $r > 0$ and $x' \in E$ with $\emph{dist}(x',X) \ge r$ define
\[C_n = \overline{\emph{co}}\left(X \cup \{x'\}\right) \setminus \bigcup_{x \in X}B(x,\emph{dist}(x',X) - 1/n).\]
Then, the sequence $\left(\emph{diam}C_n\right)_{n \in \mathbb{N}}$
converges to $0$ uniformly with respect to $x' \in E$ such that $\emph{dist}(x',X) \ge r$.
\end{lemma}
\begin{proof} 
Let $r > 0$ and $\epsilon > 0$. Take $n_0 \in \mathbb{N}$ such that
\[\frac{1}{n_0} < \frac{r\epsilon}{4(r + \mbox{diam}X)}.\]
Let $x' \in E$ with $\mbox{dist}(x',X) \ge r$. We show that for $n \ge n_0 $, $C_{n} \subseteq B(x', \epsilon/2)$. For $y \in C_{n}$ there exists $x_y \in X$ such that
\[d(y,x_y) - \frac{1}{n} < \mbox{dist}(y,X).\]
Suppose $y \notin B(x',\epsilon/2)$. Applying Lemma \ref{Lemma_ineq_co},
\begin{align*}
\mbox{dist}(y,X) & \le \mbox{dist}(x',X) - \frac{\mbox{dist}(x',X)}{\mbox{dist}(x',X) + \mbox{diam}X}d(x',y) \\
& \le \mbox{dist}(x',X) - \frac{r}{r + \mbox{diam}X}d(x',y) < \mbox{dist}(x',X) - \frac{2}{n}.
\end{align*}
Hence, 
\[d(y,x_y) < \mbox{dist}(x',X) - 1/n\]
and so $y \in B\left(x_y,\mbox{dist}(x',X) - 1/n\right)$ which is false. This means 
\[y \in B(x',\epsilon/2) \mbox{ and } C_{n} \subseteq B(x',\epsilon/2).\] 
Therefore, it follows that diam$C_n < \epsilon,$ which completes the proof.
\end{proof}

In order to state our main results, we introduce the following notations. Let $A \in P_{b,cl}(E)$ be fixed. Then we can denote $\lambda_X = \lambda_{XA}$ for $X \in P_{b,cl}(E)$. Following \cite{DeBlasi}, set
\[P_{b,cl,cv}^A(E) = \overline{\left\{X \in P_{b,cl,cv}(E) : \lambda_{X} > 0\right\}}.\]
Together with the Pompeiu-Hausdorff distance, $P_{b,cl,cv}^A(E)$ is a complete metric space if the metric of $E$ is convex.\\
For $p \in \mathbb{N}$ define
\[\mathcal{L}_p = \left\{X \in P_{b,cl,cv}^A(E) : \inf_{\sigma > 0}\mbox{diam}L_{XA}(\sigma) < \frac{1}{p} \mbox{ and } \inf_{\sigma > 0}\mbox{diam}L_{AX}(\sigma) < \frac{1}{p} \right\}\]
and
\[\mathcal{M}_p = \left\{X \in P_{b,cl,cv}(E) : \inf_{\sigma > 0}\mbox{diam}M_{XA}(\sigma) < \frac{1}{p} \mbox{ and } \inf_{\sigma > 0}\mbox{diam}M_{AX}(\sigma) < \frac{1}{p} \right\}.\]

We prove next the two main results of this section which are counterparts in the geodesic case of \cite[Theorem 3.3]{DeBlasi} and \cite[Theorem 4.3]{DeBlasi} respectively.
\begin{theorem} \label{Th_min_b+cl}
Let $E$ be a complete geodesic metric space with convex metric and curvature bounded below by $\kappa < 0$. Suppose $A \in P_{b,cl}(E)$. Then the set
\[\mathcal{W}_{min} = \left\{ X \in P_{b,cl,cv}^A(E) : \min(A,X) \mbox{ is well-posed}\right\}\]
is a dense $G_\delta$-set in $P_{b,cl,cv}^A(E)$.
\end{theorem}
\begin{proof}
Applying Proposition \ref{Prop_well-posed}, it is immediate that 
\[\mathcal{W}_{min} = \bigcap_{p \in \mathbb{N}}\mathcal{L}_p.\]
Hence, the conclusion follows if we prove that for every $p \in \mathbb{N}, \mathcal{L}_p$ is dense and open in $P_{b,cl,cv}^A(E)$. Let $p \in \mathbb{N}$. 

We first show that $\mathcal{L}_p$ is dense in $P_{b,cl,cv}^A(E)$. Take $X \in P_{b,cl,cv}^A(E)$ and $r > 0$. We want to prove that there exists $Y \in \mathcal{L}_p$ such that $h(X,Y) \le r$. Obviously, we can suppose that $r < \lambda_X$.\\
By Proposition \ref{Prop_F_kappa}, it follows that there exists $\sigma < r$ such that $F_\kappa(r, \lambda_X, \sigma) + 2\sigma < 1/p$. This clearly implies that $\sigma < 1/p$.\\
Using Lemma $\ref{Lemma_conv_diam}$ we obtain that there exists $n_0 \in \mathbb{N}$ such that $\mbox{diam}C_n < \sigma/2$ for each $n \ge n_0$ and $x' \in E$ with $\mbox{dist}(x',X) \ge r/2$, where $C_n$ is as in Lemma \ref{Lemma_conv_diam}.\\
Let $\tau = \min\{\sigma/2, 1/n_0\}$ and let $x_1 \in X$ and $a_1 \in A$ be such that 
\[d(x_1,a_1) < \lambda_X + \frac{\tau}{2}.\]
Take $x' \in [x_1,a_1]$ such that $d(x_1,x') = r$. Now consider $Y = \overline{\mbox{co}}\left(X \cup \{x'\}\right)$. Then it is easy to see that $h(X,Y) \le r$.\\ 
We also have that 
\[\lambda_Y \le \mbox{dist}(x',A) \le d(x',a_1) = d(x_1,a_1) - r < \lambda_X + \frac{\tau}{2} - r.\]
Likewise,
\[\mbox{dist}(x',X) \ge \lambda_X - \mbox{dist}(x',A) \ge r - \frac{\tau}{2} \ge \frac{r}{2}.\]
We show next that $\lambda_Y \ge \lambda_X - r$. This would yield $Y \in P_{b,cl,cv}^A(E)$. Suppose there exist $y \in Y$ and $a \in A$ such that $d(y,a) < \lambda_X - r$. By Lemma \ref{Lemma_ineq_co},
\[\mbox{dist}(y,X) \le r - 2\alpha d(x',y) \mbox{ where } \alpha = \frac{r/2}{2\left(r/2 + \mbox{diam}X\right)}.\]
Suppose $y \ne x'$. Then there is $x^* \in X$ such that $d(y,x^*) < \mbox{dist}(y,X) + \alpha d(x',y).$ Thus,
\begin{align*}
\lambda_X & \le d(x^*,a) \le d(x^*,y) + d(y,a) < \mbox{dist}(y,X) + \alpha d(x',y) + \lambda_X - r \\
& \le r - \alpha d(x',y) + \lambda_X - r < \lambda_X.
\end{align*}
This is false, so $y = x'.$ In this case, $d(x',a) < \lambda_X - r$ implies $d(x_1,a) < \lambda_X$ which is a contradiction. Therefore, $\lambda_Y \ge \lambda_X - r$.\\
Let $y \in L_{YA}(\tau/2) \setminus C_{n_0}$. Then,
\[\mbox{dist}(y,A) \le \frac{\tau}{2} + \lambda_Y \mbox{ and } \mbox{dist}(y,X) < \mbox{dist}(x',X) - \frac{1}{n_0}.\]
This implies
\begin{align*}
\lambda_X & \le \mbox{dist}(y,A) + \mbox{dist}(y,X) < \frac{\tau}{2} + \lambda_Y + \mbox{dist}(x',X) - \frac{1}{n_0} \le \frac{\tau}{2} + \lambda_Y + r - \frac{1}{n_0}\\
&< \lambda_X + \tau -\frac{1}{n_0} \le \lambda_X
\end{align*}
But this is a contradiction, so $L_{YA}(\tau/2) \subseteq C_{n_0}$ which means
 \[\mbox{diam}L_{YA}\left(\frac{\tau}{2}\right) < \frac{\sigma}{2} < \frac{1}{p}.\]
Let $a \in L_{AY}(\tau/4)$. Then $\mbox{dist}(a,Y) \le \lambda_Y + \tau/4$.\\
Pick $y \in Y$ such that $d(a,y) \le \lambda_Y + \frac{\tau}{2} < \lambda_X + \tau -r$. This yields
\[\mbox{dist}(y,X) \ge \lambda_X - \mbox{dist}(y,A) > r - \tau \ge \mbox{dist}(x',X) - \frac{1}{n_0}.\]
Consequently, $y \in C_{n_0}$ and $d(y,x') < \sigma/2$. This means
\[d(a,x') \le d(a,y) + d(y,x') < \lambda_X + \tau - r + \frac{\sigma}{2} \le \lambda_X - r + \sigma,\]
and so $a \in \widetilde{B}(x',\lambda_X-r+\sigma)$. Since $a \notin B(x_1,\lambda_X)$ it is clear that 
\[a \in D(x_1,x';\lambda_X,\sigma).\]
Applying Proposition \ref{Prop_diam}, we obtain that 
\[\mbox{diam}L_{AY}\left(\frac{\tau}{4}\right) \le F_\kappa(r,\lambda_X,\sigma) + 2\sigma < \frac{1}{p}.\]
This completes the proof that $Y \in \mathcal{L}_p$.

Let us now show that $\mathcal{L}_p$ is open. Consider $X \in \mathcal{L}_p$ and let 
\[\theta = \max\left\{\inf_{\sigma > 0}\mbox{diam}L_{XA}(\sigma), \inf_{\sigma > 0}\mbox{diam}L_{AX}(\sigma)\right\}.\]
Choose $\epsilon > 0$ such that $\theta + 2\epsilon < 1/p$. Also, let $\sigma > 0$ be such that
\[\max\left\{\mbox{diam}L_{XA}(\sigma), \mbox{diam}L_{AX}(\sigma)\right\} < \theta + \epsilon.\]
Take $\delta = \min\{\sigma/4,\epsilon/2\}$. Let $Y \in P_{b,cl,cv}^A(E)$ with $h(X,Y) < \delta$. We show that $Y \in \mathcal{L}_p$. Since $h(X,Y) < \delta$ it clearly follows that $\lambda_Y < \lambda_X + \delta$. Let $y \in L_{YA}(\sigma/2)$. Then there exists $x \in X$ such that $d(x,y) < \delta$. We also have that $x \in L_{XA}(\sigma)$ because
\[\mbox{dist}(x,A) \le \mbox{dist}(y,A) + h(X,Y) < \lambda_Y + \frac{\sigma}{2} + \delta < \lambda_X + 2\delta + \frac{\sigma}{2} \le \lambda_X + \sigma.\]
For $y_1, y_2 \in L_{YA}(\sigma/2)$ arbitrary, there exist $x_1, x_2 \in L_{XA}(\sigma)$ such that $d(x_1,y_1) < \delta$ and $d(x_2,y_2) < \delta$. Thus, $\mbox{diam}L_{AY}(\sigma/2) < 1/p$  since
\[d(y_1,y_2) \le 2\delta + \mbox{diam}L_{XA}(\sigma) < \theta + 2\epsilon < \frac{1}{p}.\]
Let $a \in L_{AY}(\sigma/2)$. Then
\[\mbox{dist}(a,X) \le \mbox{dist}(a,Y) + h(X,Y) < \lambda_Y + \frac{\sigma}{2} + \delta < \lambda_X + 2\delta + \frac{\sigma}{2} \le \lambda_X + \sigma.\]
This yields $a \in L_{AX}(\sigma)$ and so $\mbox{diam}L_{AY}(\sigma/2) < \theta + \epsilon < 1/p.$\\
Hence, $Y \in \mathcal{L}_p$ and the proof is complete.
\end{proof}

In the sequel we give the corresponding maximization result.

\begin{theorem} \label{Th_max_b+cl}
Let $E$ be a complete geodesic metric space with convex metric, the geodesic extension property and curvature bounded below by $\kappa < 0$. Suppose $A \in P_{b,cl}(E)$. Then the set
\[\mathcal{W}_{max} = \left\{ X \in P_{b,cl,cv}(E) : \max(A,X) \mbox{ is well-posed}\right\}\]
is a dense $G_\delta$-set in $P_{b,cl,cv}(E)$.
\end{theorem}
\begin{proof}
By Proposition \ref{Prop_well-posed}, it is immediate that 
\[\mathcal{W}_{max} = \bigcap_{p \in \mathbb{N}}\mathcal{M}_p.\]
Again, the conclusion follows if we prove that for every $p \in \mathbb{N}, \mathcal{M}_p$ is dense and open in $P_{b,cl,cv}(E)$. Let $p \in \mathbb{N}$. 

We first show that $\mathcal{M}_p$ is dense in $P_{b,cl,cv}(E)$. Take $X \in P_{b,cl,cv}(E)$ and $r > 0$. We want to prove that there exists $Y \in \mathcal{M}_p$ such that $h(X,Y) \le r$. Obviously, we can suppose that $r < \mu_X$.\\
By Proposition \ref{Prop_F_kappa}, it follows that there exists $\sigma < r$ such that $F_\kappa(r, \mu_X + r - \sigma, \sigma) + 2\sigma < 1/p$. This clearly implies that $\sigma < 1/p$.\\
Using Lemma $\ref{Lemma_conv_diam}$ we obtain that there exists $n_0 \in \mathbb{N}$ such that $\mbox{diam}C_n < \sigma/2$ for each $n \ge n_0$ and $x' \in E$ with $\mbox{dist}(x',X) \ge r/2$, where $C_n$ is as in Lemma \ref{Lemma_conv_diam}.\\
Let $\tau = \min\{\sigma/2, 1/n_0\}$ and let $x_1 \in X$ and $a_1 \in A$ be such that 
\[d(x_1,a_1) > \mu_X - \frac{\tau}{4} > \mu_X - r > 0.\]
Since $E$ has the geodesic extension property  there exists a point $x'$ on the geodesic line determined by $x_1$ and $a_1$ such that $d(x_1,x') = r$ and $d(a_1,x') = r + d(x_1,a_1)$. Now consider $Y = \overline{\mbox{co}}\left(X \cup \{x'\}\right)$. It is easy to see that $h(X,Y) \le r$.\\ 
The following holds 
\[\mu_Y \ge \mbox{Dist}(x',A) \ge d(x',a_1) = d(x_1,a_1) + r > \mu_X - \frac{\tau}{4} + r.\]
Likewise,
\[\mbox{dist}(x',X) \ge  \mbox{Dist}(x',A) - \mu_X \ge r - \frac{\tau}{4} \ge \frac{r}{2}.\]
Let $y \in M_{YA}(\tau/2) \setminus C_{n_0}$. Then,
\[\mbox{Dist}(y,A) \ge \mu_Y - \frac{\tau}{2} \mbox{ and } \mbox{dist}(y,X) < \mbox{dist}(x',X) - \frac{1}{n_0}.\]
These inequalities imply
\[\mbox{Dist}(y,A) \ge \mu_Y- \frac{\tau}{2} > \mu_X + r - \frac{3}{4}\tau \ge \mu_X + r  - \frac{3}{4n_0},\]
and 
\[\mbox{Dist}(y,A) \le \mu_X + \mbox{dist}(y,X) < \mu_X + \mbox{dist}(x',X) - \frac{1}{n_0} \le \mu_X + r - \frac{1}{n_0},\]
which taken together yield a contradiction. Therefore, $M_{YA}(\tau/2) \subseteq C_{n_0}$ which means
 \[\mbox{diam}M_{YA}(\tau/2) < \frac{\sigma}{2} < \frac{1}{p}.\]
Let $a \in M_{AY}\left(\frac{\tau}{4}\right)$. Then $\mbox{Dist}(a,Y) \ge \mu_Y - \tau/4$.\\
Pick $y \in Y$ such that 
\[d(a,y) \ge \mu_Y - \tau/2 > \mu_X + r - \frac{3}{4}\tau.\] 
This yields
\[\mbox{dist}(y,X) \ge \mbox{Dist}(y,A) - \mu_X > r - \frac{3}{4}\tau > \mbox{dist}(x',X) - \frac{1}{n_0}.\]
Consequently, $y \in C_{n_0}$ and $d(y,x') < \sigma/2$. This means
\[d(a,x') \ge d(a,y) - d(y,x') > \mu_X + r - \frac{3}{4}\tau - \frac{\sigma}{2} > \mu_X + r - \sigma,\]
and so $a \notin B(x',\mu_X + r - \sigma)$. Since $a \in \widetilde{B}(x_1,\mu_X)$ it is clear that 
\[a \in D(x',x_1;\mu_X + r - \sigma,\sigma).\] 
Applying Proposition \ref{Prop_diam}, we obtain that 
\[\mbox{diam}M_{AY}\left(\frac{\tau}{4}\right) \le F_\kappa(r,\mu_X + r - \sigma,\sigma) + 2\sigma < \frac{1}{p}.\]
This completes the proof that $Y \in \mathcal{M}_p$.

The fact that $\mathcal{M}_p$ is open follows in a similar manner as in the proof of Theorem \ref{Th_min_b+cl}.
\end{proof}

We conclude this section by giving a characterization of the well-posedness of the problem $\min(X,Y)$ in complete CAT$(0)$ spaces. We prove that in the following particular context, the conditions in Proposition \ref{Prop_well-posed} can be relaxed.

\begin{proposition}
Let $E$ be a complete CAT$(0)$ space, $X \in P_{b,cl,cv}(E)$ and $Y \in P_{b,cl}(E)$. The problem $\min(X,Y)$ is well-posed if and only if
\[\inf_{\sigma > 0}\emph{diam}L_{YX}(\sigma) = 0.\]
\end{proposition}
\begin{proof}
We solely need to prove the sufficiency part. Suppose $\inf_{\sigma > 0}\mbox{diam}L_{YX}(\sigma) = 0$. Take $(x_n,y_n)_{n \in \mathbb{N}}$ in $X \times Y$ a minimizing sequence. Then, for every $\epsilon > 0$ there exists $n_0$ such that $y_n \in L_{YX}(\epsilon)$ for every $n \ge n_0$. This implies that the sequence $(y_n)_{n \in \mathbb{N}}$ is Cauchy and hence it converges to some $y \in Y$ for which $\lim_{n \to \infty}d(x_n,y) = \lambda_{XY}$. Since $(x_n)_{n \in \mathbb{N}}$ is a sequence in a bounded, closed and convex subset of a complete CAT$(0)$ space, it contains a $\Delta$-convergent subsequence, $(x_{n_k})_{k \in \mathbb{N}}$, whose $\Delta$-limit, $x$, belongs to $X$. The generalized (CN) inequality yields that, for all $t \in [0,1]$, we have that
\[d(x_{n_k},(1-t)x+ty)^2 \le (1-t)d(x_{n_k},x)^2 + td(x_{n_k},y)^2 - t(1-t)d(x,y)^2.\]
Taking the superior limit with respect to $k$ and using the fact that 
\[\limsup_{k \to \infty}d(x_{n_k},x) \le \limsup_{k \to \infty}d(x_{n_k}, (1-t)x + ty),\]
we obtain that
\[\limsup_{k \to \infty}d(x,x_{n_k})^2 \le \limsup_{k \to \infty}d(x_{n_k},y)^2 - (1-t)d(x,y)^2.\]
Letting $t \searrow 0$ implies that
\[d(x,y)^2 + \limsup_{k \to \infty}d(x_{n_k},x)^2 \le \lambda_{XY}^2.\]
Hence, $d(x,y) = \lambda_{XY}$ and $\lim_{k \to \infty}x_{n_k} = x.$ Thus, $x \in P_X(y)$. The set $P_X(y)$ is a singleton in this setting (for a justification see for instance \cite[Proposition 2.4]{Bridson}) and so $\{x\} = P_X(y)$. This means that the sequence $(x_n)_{n \in \mathbb{N}}$ is $\Delta$-convergent to $x$ and, as before, one can see that $\lim_{n \to \infty}x_n = x.$\\
To prove the uniqueness of the solution, suppose there exists $(x_0,y_0)$ another solution. Then the sequence $(x,y), (x_0,y_0), (x,y), (x_0,y_0), \ldots$ is a minimizing sequence. However, based on the above, the sequences $x,x_0,x,x_0,\ldots$ and $y,y_0,y,y_0,\ldots$ must be convergent, so $x = x_0$ and $y = y_0$. This completes the proof that $\min(X,Y)$ is well-posed.
\end{proof}

\section{Results involving compactness}
In this section we study the same problems as in section $3$ but we modify conditions we imposed in our results. More particularly, we focus on the situation in which the set $A$ is compact. We show that under this stronger assumption on the set we can weaken the condition on the geodesic space from being of curvature bounded below to not having bifurcating geodesics. However, in the first theorem we need to add the reflexivity condition on the space. Before stating this result we give the following property whose proof follows similar patterns as in \cite[Proposition 7]{Dhomp&Kirk}.

\begin{lemma} \label{Lemma_refl_compl}
Let $(E,d)$ be a reflexive geodesic metric space with convex metric. Then $E$ is complete.
\end{lemma}
\begin{proof}
Let $(x_n)_{n \in \mathbb{N}}$ be a Cauchy sequence, define $\varphi : E \to \mathbb{R_+}$, $\varphi(u) = \limsup_{n \to \infty}d(u,x_n)$ and set $r = \inf_{u \in E}\varphi(u)$.\\
Let $p \in \mathbb{N}$. Then there exists $u_p \in E$ such that $\varphi(u_p) < r + 1/p$. Hence, for $n$ sufficiently large, $d(u_p,x_n) < r + 1/p$, that is, $u_p \in B(x_n,r + 1/p)$. Thus,
\[C_p = \bigcup_{k \in \mathbb{N}} \left(\bigcap_{i \ge k}B(x_i,r + 1/p)\right) \ne \emptyset.\]
The set $C_p$ is bounded since the sequence $(x_n)_{n \in \mathbb{N}}$ is bounded. Because the metric of $E$ is convex, the sets $C_p$ and $\overline{C_p}$ will be convex. Hence, the sequence $(\overline{C_p})_{p \in \mathbb{N}}$ is a descending sequence of nonempty, bounded, closed and convex sets and so, by the reflexivity of the space, there exists $u \in E$ such that
\[u \in \bigcap_{p \in \mathbb{N}}C_p.\]
This yields that $\varphi(u) \le r$ which means $\limsup_{n \to \infty}d(u,x_n) \le \limsup_{n \to \infty}d(x_m,x_n)$ for every $m \in \mathbb{N}$. But since the sequence $(x_n)_{n \in \mathbb{N}}$ is Cauchy this implies that $\lim_{n \to \infty}d(u,x_n) = 0$. Therefore, the sequence $(x_n)_{n \in \mathbb{N}}$ is convergent and the proof ends.
\end{proof}

\begin{theorem} \label{Th_min_cp}
Let $E$ be a reflexive geodesic space with convex metric and no bifurcating geodesics. Suppose $A \in P_{cp}(E)$. Then the set
\[\mathcal{W}_{min} = \left\{ X \in P_{b,cl,cv}^A(E) : \min(A,X) \mbox{ is well-posed}\right\}\]
is a dense $G_\delta$-set in $P_{b,cl,cv}^A(E)$.
\end{theorem}
\begin{proof}
We begin by proving the denseness result. Take $X \in P_{b,cl,cv}^A(E)$. First, let us show that there exist $a_1 \in A$ and $x_1 \in X$ such that $d(a_1,x_1) = \lambda_X$. Take $(a'_n,x'_n)_{n \in \mathbb{N}}$ in $A \times X$ such that $\lim_{n \to \infty}d(a'_n,x'_n) = \lambda_X$. Since $A$ is compact, $(a'_n)_{n \in \mathbb{N}}$ contains a convergent subsequence, $(a'_{n_k})_{k \in \mathbb{N}}$. Suppose $\lim_{k \to \infty}a'_{n_k} = a_1 \in A.$ Then $\lim_{k \to \infty}d(a_1,x'_{n_k}) = \lambda_X$. Now let
\[C_k = \left\{x \in X : d(a_1,x) \le \lambda_X + \frac{1}{k}\right\}.\]
It is clear that $(C_k)_{k \in \mathbb{N}}$ is a descending sequence of nonempty, bounded and closed sets. The convexity of the metric assures that these sets are also convex. Using the reflexivity of the space, we obtain that there exists $x_1 \in X$, 
\[x_1 \in \bigcap_{k \in \mathbb{N}}C_k.\]
This yields that $d(a_1,x_1) = \lambda_X$.\\
Let $r > 0$. We want to prove that there exists $Y \in P_{b,cl,cv}^A(E)$ such that $h(X,Y) \le r$ and $\min(A,Y)$ is well-posed. Obviously, we can suppose that $r < \lambda_X$.\\ 
Take $x' \in [x_1,a_1]$ such that $d(x_1,x') = r$. Then, as it was shown in \cite{Zamfirescu} by Zamfirescu, the fact that the space has no bifurcating geodesics guarantees that $P_A(x') = \{a_1\}$. We write next the proof of this for the sake of completeness. Suppose $a \in P_A(x').$ Then 
\[d(x_1,a_1) \le d(x_1,a) \le d(x_1,x') + d(x',a) \le d(x_1,x') + d(x',a_1) = d(x_1,a_1),\]
and so $d(x_1,a) = d(x_1,x') + d(x',a)$, which means $x' \in [x_1,a]$. But if $a \ne a_1$ then this will contradict the fact that the space has no bifurcating geodesics.\\
Now consider $Y = \overline{\mbox{co}}\left(X \cup \{x'\}\right)$. Then it is easy to see that $h(X,Y) \le r$. It is also clear that $\lambda_Y \le \lambda_X - r$ since $d(x',a_1) = \lambda_X - r$.\\
Let $(a_n,y_n)_{n \in \mathbb{N}}$ be a sequence in $A \times Y$ such that $\lim_{n \to \infty}d(a_n,y_n) = \lambda_Y$. Suppose $\beta = \limsup_{n \to \infty}d(x',y_{n})$ and \[\alpha = \frac{\mbox{dist}(x',X)}{\mbox{dist}(x',X) + \mbox{diam}X} > 0.\]
Applying Lemma \ref{Lemma_ineq_co},
\[\lambda_X \le \mbox{dist}(a_{n}, X) \le d(a_{n}, y_{n}) + \mbox{dist}(y_{n},X) \le d(a_{n}, y_{n}) + \mbox{dist}(x',X) - \alpha d(y_{n},x').\]
This implies that $\alpha \beta \le \lambda_Y + r - \lambda_X \le 0$ and so $\lambda_Y = \lambda_X - r > 0$ and $\beta = 0$. 
This means on the one hand that $Y \in P_{b,cl,cv}^A(E)$ and on the other that $\lim_{n \to \infty}y_{n} = x'$.\\
 Because $A$ is compact, the sequence $(a_n)_{n \in \mathbb{N}}$ has a convergent subsequence $(a_{n_k})_{k \in \mathbb{N}}$. Suppose $\lim_{k \to \infty}a_{n_k} = a$ for some $a \in A.$ Then it follows that $d(x',a) = \lambda_Y$ and hence $a \in P_A(x') = \{a_1\}$, that is, $a = a_1$. Based on the compactness of $A$, we conclude that the sequence $(a_n)_{n \in \mathbb{N}}$ converges to $a_1$.\\
Suppose there exists another solution, say $(a,y)$, of the problem $\min(A,Y)$. Then the sequence $(a,y), (a_1,x'), (a,y), (a_1,x'), \ldots$ is a minimizing sequence. However, based on the above, the sequences $a,a_1,a,a_1,\ldots$ and $y,x',y,x',\ldots$ must be convergent, so $a = a_1$ and $y = x'$. This completes the proof that $\min(A,Y)$ is well-posed. Notice that this implies that $Y \in \mathcal{L}_p$ for each $p \in \mathbb{N}$.

The fact that $\mathcal{W}_{min}$ is a $G_\delta$ set can be easily seen by applying Proposition \ref{Prop_well-posed} (here we use the completeness of the space which is assured by Lemma \ref{Lemma_refl_compl}) and writing 
\[\mathcal{W}_{min} = \bigcap_{p \in \mathbb{N}}\mathcal{L}_p.\]
Then we prove similarly as in Theorem \ref{Th_min_b+cl} that the set $\mathcal{L}_p$ is open in $P_{b,cl,cv}^A(E)$ for every $p \in \mathbb{N}$.
\end{proof}

The following is a particular case of the above result.
\begin{corollary}
Let $E$ be a complete CAT$(0)$ space with no bifurcating geodesics. Suppose $A \in P_{cp}(E)$. Then the set
\[\mathcal{W}_{min} = \left\{ X \in P_{b,cl,cv}^A(E) : \min(A,X) \mbox{ is well-posed}\right\}\]
is a dense $G_\delta$-set in $P_{b,cl,cv}^A(E)$.
\end{corollary}

\begin{remark}
The proof of Theorem \ref{Th_min_cp} relies on the fact that $\min(A,X)$ always has a solution. In fact, the reflexivity of the space is mainly used to ensure this condition. Can we drop the condition that the problem has a solution?
\end{remark}

Next we focus on the maximization problem for $A$ compact. In order to follow the same line of argument as in the previous result we need the fact that the problem $\max(A,X)$ has a solution. However, in \cite{Sababheh} it is proved that in a reflexive Banach space, the remotal distance from a point to a bounded, closed and convex set is guaranteed to be reached if and only if the space is finite dimensional. This is why we impose the compactness condition on the set $X$ in our next result in order to make sure that $\max(A,X)$ finds a solution.

\begin{theorem}
Let $E$ be a complete geodesic space with no bifurcating geodesics and the geodesic extension property. Suppose $A \in P_{cp}(E)$. Then the set
\[\mathcal{W}_{max} = \left\{ X \in P_{cp}(E) : \max(A,X) \mbox{ is well-posed}\right\}\]
is a dense $G_\delta$-set in $P_{cp}(E)$.
\end{theorem}
\begin{proof}
We will briefly sketch the proof of the denseness result. Take $X \in P_{cp}(E)$ and $r > 0$. We want to prove that there exists $Y \in P_{cp}(E)$ such that $h(X,Y) \le r$. Obviously, we can suppose that $r < \mu_X$.\\
Since the sets $A$ and $X$ are both compact, there exist $a_1 \in A$ and $x_1 \in X$ such that $d(a_1,x_1) = \mu_X$.\\
Because $E$ has the geodesic extension property there exists a point $x'$ on the geodesic line determined by $x_1$ and $a_1$ such that $d(x_1,x') = r$ and $d(a_1,x') = r + d(x_1,a_1)$. Then the fact that the space has no bifurcating geodesics guarantees that $F_A(x') = \{a_1\}$. Indeed, suppose $a \in F_A(x').$ Then 
\[d(x',a_1) = d(x',x_1) + d(x_1,a_1)\ge d(x',x_1) + d(x_1,a) \ge d(x',a) \ge d(x',a_1),\]
and so $d(x',a) = d(x',x_1) + d(x_1,a)$, which means $x_1 \in [x',a]$. But if $a \ne a_1$ then this will contradict the fact that the space has no bifurcating geodesics.\\
Now consider $Y = X \cup \{x'\}$. Clearly, $h(X,Y) \le r$, $Y \in P_{cp}(E)$ and it is easy to prove that $\max(A,Y)$ is well-posed.
\end{proof}

\begin{remark} 
In the above result we do not need the completeness of the space to prove the denseness result because the argument avoids the use of Proposition \ref{Prop_well-posed} and of the Baire category theorem. Thus, this proof follows a more direct approach than the ones stated in section \ref{Results_curv_b_b}. 
\end{remark}

\begin{remark}
Regarding the problem $\max(A,X)$, where the fixed set $A$ is compact, we raise the following question: is the set  
\[\mathcal{W}_{max} = \left\{ X \in P_{cp,cv}(E) : \max(A,X) \mbox{ is well-posed}\right\}\]
a dense $G_\delta$-set in $P_{cp,cv}(E)$? The Hopf-Rinow Theorem (see \cite[Proposition 3.7]{Bridson}) states that if $E$ is complete and locally compact, then it is proper. Hence, if the space is additionally locally compact and with convex metric then we can answer the question in the positive by taking in the above proof the set $Y= \overline{\emph{co}}\left(X \cup \{x'\}\right)$, which will be a compact and convex set.
\end{remark}

\section{The Drop Theorem in spaces with convex metric}
In \cite{Danes}, Dane\v s proved the following geometric result known as the Drop Theorem.
\begin{theorem}
Let $(E,\|$ $\|)$ be a Banach space. Suppose $B$ is the unit ball in $E$ and let $A \in P_{cl}(E)$ be such that $\inf\{\| x \| : x \in A\} > 1.$ Then there exists $a \in A$ such that 
\[\emph{co}\left(B \cup \{a\}\right) \cap A = \{a\},\]
where $\emph{co}\left(B \cup \{a\}\right)$ denotes the convex hull in the Banach space of the set $B \cup \{a\}$.
\end{theorem}
The name of this theorem has its origin from the fact that the set $\mbox{co}\left(B \cup \{a\}\right)$ was called a drop. Equivalences of this result or of its generalized versions with other fundamental theorems in nonlinear analysis and various areas of their applications are discussed for instance in \cite{Georgiev, Penot}.

In this section we prove a variant of the Drop Theorem in the setting of a geodesic space with convex metric. We will derive this result from the following theorem called the Strong Flower Petal Theorem. For a proof of this theorem see \cite[Proposition 2.5]{Georgiev}. This result uses the following extension of the definition of a petal given in \cite{Penot}: for $(E,d)$ a metric space and a function $f : E \to \mathbb{R}$, we say that the set
\[P_{\alpha,\delta}(x_0,f) = \left\{x \in E : f(x) \le f(x_0) - \alpha d(x,x_0) + \delta\right\}\]
is the petal associated to $\delta \ge 0$, $\alpha > 0$, $x_0 \in E$ and $f$.

\begin{theorem} [Strong Flower Petal Theorem]\label{Th_strongPetal}
Let $(E,d)$ be a complete metric space, $A \in P_{cl}(E)$ and $f : E \to \mathbb{R}$ a Lipschitz function bounded below on $A$. Suppose $\delta > 0$, $\alpha > 0$ and $x_0 \in A$. Then there exists a point $a \in A \cap P_{\alpha,\delta}(x_0,f)$ such that
\begin{itemize}
\item[\emph{(i)}] $P_{\alpha,0}(a,f) \cap A = \{a\}$ and
\item[\emph{(ii)}] $x_n \to a$ for every sequence $(x_n)_{n \in \mathbb{N}}$ in $P_{\alpha,0}(a,f)$ with $\emph{dist}(x_n,A) \to 0.$
\end{itemize}
\end{theorem}

The following is a variant of the Drop Theorem in geodesic spaces.

\begin{theorem} \label{Th_drop}
Let $(E,d)$ be a complete geodesic space with convex metric and let $A \in P_{cl}(E)$ and $B \in P_{b,cl,cv}(E)$ be such that $\lambda_{AB} > 0$. Suppose $\epsilon > 0$. Then there exists $a \in A$ such that
\begin{itemize}
\item[\emph{(i)}] $\emph{dist}(a,B) < \lambda_{AB} + \epsilon$,
\item[\emph{(ii)}] $\overline{\emph{co}}\left(B \cup \{a\}\right) \cap A = \{a\}$ and
\item[\emph{(iii)}] $x_n \to a$ for every sequence $(x_n)_{n \in \mathbb{N}}$ in $\overline{\emph{co}}\left(B \cup \{a\}\right)$ with $\emph{dist}(x_n,A) \to 0.$
\end{itemize}
\end{theorem}
\begin{proof}
Let $\epsilon > 0$. Then there exists $x_0 \in A$ such that $\mbox{dist}(x_0,B) < \lambda_{AB} + \epsilon/2$. Take
\[\alpha = \frac{\lambda_{AB}}{\lambda_{AB} + \mbox{diam}B} > 0.\]
Since the function $\mbox{dist}(\cdot,B)$ is nonexpansive we can apply Theorem \ref{Th_strongPetal} which yields that there exists $a \in A \cap P_{\alpha,\epsilon/2}\left(x_0,\mbox{dist}(\cdot,B)\right)$ such that
\begin{itemize}
\item[\mbox{(a)}] $P_{\alpha,0}\left(a,\mbox{dist}(\cdot,B)\right) \cap A = \{a\}$ and
\item[\mbox{(b)}] $x_n \to a$ for every sequence $(x_n)_{n \in \mathbb{N}}$ in $P_{\alpha,0}\left(a,\mbox{dist}(\cdot,B)\right)$ with $\mbox{dist}(x_n,A) \to 0.$
\end{itemize}
Since $a \in P_{\alpha,\epsilon/2}\left(x_0,\mbox{dist}(\cdot,B)\right)$ it follows that $\mbox{dist}(a,B) < \lambda_{AB} +\epsilon$ and thus (i) holds.\\
Let us prove that 
\begin{equation} \label{incl_drop}
\overline{\mbox{co}}\left(B \cup \{a\}\right) \subseteq P_{\alpha,0}\left(a,\mbox{dist}(\cdot,B)\right).
\end{equation}
Take $y \in \overline{\mbox{co}}\left(B \cup \{a\}\right)$. Applying Lemma \ref{Lemma_ineq_co} we obtain that,
\[\mbox{dist}(y,B) \le \mbox{dist}(a,B) - \frac{\mbox{dist}(a,B)}{\mbox{dist}(a,B) + \mbox{diam}B}d(a,y) \le \mbox{dist}(a,B) - \alpha d(a,y),\]
and hence $y \in P_{\alpha,0}\left(a,\mbox{dist}(\cdot,B)\right)$. Now it is immediate that conditions (a) and (b) imply (ii) and (iii) respectively.
\end{proof}

As an application of this version of the Drop Theorem we will obtain an analogue of an optimization result proved by Georgiev \cite[Theorem 4.2]{Georgiev} in the context of Banach spaces. In order to state this result we need to briefly introduce some notions which can also be found in \cite{Georgiev}.

Let $(E,d)$ be complete metric space, $f:E \to \mathbb{R}$ a lower semi-continuous and bounded below function and $A \in P_{b,cl}(E)$.
The minimization problem denoted by $\min(A,f)$ consists in finding $x_0 \in A$ (the solution of the problem) such that $f(x_0) = \inf\{f(x) : x \in A\}$.\\
For $\sigma > 0$, let 
\[L_{A,f}(\sigma) = \left\{x \in E : f(x) \le \inf_{y \in A}f(y) + \sigma \mbox{ and } \mbox{dist}(x,A) \le \sigma \right\}.\]
The problem $\min(A,f)$ is well-posed in the sense of Revalski if $\inf_{\sigma > 0}\mbox{diam}L_{A,f}(\sigma) = 0$ (see \cite{Revalski}). This is equivalent to requesting that it has a unique solution $x_0 \in A$ and every sequence $(x_n)_{n \in \mathbb{N}}$ in $E$ converges to $x_0$ provided $f(x_n) \to f(x_0)$ and $\mbox{dist}(x_n,A) \to 0.$

The following lemma is the counterpart of \cite[Lemma 4.1]{Georgiev} for geodesic metric spaces.

\begin{lemma} \label{Lemma_dist}
Let $E$ be a geodesic space, $X$ a bounded subset of $E$ and $f : E \to \mathbb{R}$ continuous and convex. For $c \in \mathbb{R}$, let $A = \{x \in E : f(x) \le c\}.$ Suppose there exists $z \in E$ such that $f(z) < c$. Then for every $\epsilon > 0$ there exists $\delta > 0$ such that $\emph{dist}(x,A) < \epsilon$ for each $x \in X$ with $f(x) < c + \delta.$
\end{lemma}
\begin{proof}
Suppose there exists $\epsilon > 0$ such that for every $n \in \mathbb{N}$ there exists $x_n\in X$ such that
\[f(x_n) < c + \frac{1}{n} \mbox{ and } \mbox{dist}(x_n,A) \ge \epsilon.\]
Let $t, \delta \in (0,1)$ be such that
\[ (1 - t)\mbox{Dist}(z,X) < \epsilon \mbox{ and } (1 - t)f(z) < (1-t)c - t\delta.\]
Also, let $m > 1/\delta$ and $y = (1-t)z + tx_m$. Then
\[f(y) \le (1-t)f(z) + tf(x_m) < (1-t)c - t\delta + t(c + \frac{1}{m}) < c,\]
so $y \in A$. However,
\[\epsilon \le \mbox{dist}(x_m,A) \le d(x_m,y) = (1-t)d(x_m,z) \le (1-t)\mbox{Dist}(z,X) < \epsilon,\]
which is a contradiction. This completes the proof.
\end{proof}

Before proving the optimization result we define, for $p \in \mathbb{N}$ and $E$ a geodesic space, the set
\[\mathcal{L}_p = \left\{X \in P_{b,cl,cv}(E) : \inf_{\sigma > 0}\mbox{diam}L_{X,f}(\sigma) < \frac{1}{p}\right\}.\]

\begin{theorem} \label{Th_f_wellp}
Let $E$ be a complete geodesic space with convex metric and let $f : E \to \mathbb{R}$ be continuous, convex, bounded below on bounded sets and satisfying one of the following:
\begin{itemize}
\item[\emph{(i)}] $\inf_{x \in E}f(x) = -\infty$ or
\item[\emph{(ii)}] there exists $z_0 \in E$ such that $f(z_0) = \inf_{x \in E}f(x)$ and every sequence $(x_n)_{n \in \mathbb{N}}$ in E converges to $z_0$ if $f(x_n) \to f(x_0)$.
\end{itemize}
Then the set
\[\mathcal{W}_{min} = \left\{ X \in P_{b,cl,cv}(E) : \min(X,f) \mbox{ is well-posed in the sense of Revalski}\right\}\]
is a dense $G_\delta$-set in $P_{b,cl,cv}(E)$.
\end{theorem}
\begin{proof}
Let $X \in P_{b,cl,cv}(E)$ and $r > 0$. To prove the denseness result, we show that there exists $Y \in P_{b,cl,cv}(E)$ such that $h(X,Y) \le r$ and $\min(Y,f)$ is well-posed. If hypothesis (ii) holds and $z_0$ furnished by it is such that $\mbox{dist}(z_0,X) \le r$, then take $Y = \overline{\mbox{co}}\left(X \cup \{z_0\}\right)$. It is easy to see that $h(X,Y) \le r$ and $\lim_{n \to \infty}\mbox{diam}L_{Y,f}(1/n) = 0$, so $\min(Y,f)$ is well-posed and we are done.\\
Suppose (i) holds or (ii) is accomplished with $\mbox{dist}(z_0,X) > r$. Set
\[B = \{x \in E : \mbox{dist}(x,X) \le r/2\}, \mbox{ } m = \inf\{f(x) : x \in B\}\]
and
\[A = \{x \in E : f(x) \le m\}.\]
Take $x_0 \in E$ such that $f(x_0) < m$ (in case (ii) holds with $\mbox{dist}(z_0,X) > r$ take $x_0 = z_0$). Applying Lemma \ref{Lemma_dist} we obtain that $\lambda_{AB} = 0$ and so $\lambda_{AX} \le r/2$. Suppose $\lambda_{AX} < r/2.$ Then there exists $z \in A$ and $\delta > 0$ such that $B(z,\delta) \subseteq B$. Let $t \in (0,1)$ be such that
\[(1 - t)\mbox{Dist}(x_0,B) < \delta\]
and take $y = (1-t)x_0 + tz$. Then
\[d(y,z) = (1-t)d(x_0,z) \le (1-t)\mbox{Dist}(x_0,B) < \delta.\]
Thus, $y \in B$. At the same time, $f(y) \le (1-t)f(x_0) + tf(z) < m$ which is a contradiction. This means that $\lambda_{AX} = r/2 > 0$. The set $A$ is closed because $f$ is continuous. By Theorem \ref{Th_drop}, there exists $a \in A$ such that
\begin{itemize}
\item[(a)] $\mbox{dist}(a,X) < r$,
\item[(b)] $\overline{\mbox{co}}\left(X \cup \{a\}\right) \cap A = \{a\}$ and
\item[(c)] $x_n \to a$ for every sequence $(x_n)_{n \in \mathbb{N}}$ in $\overline{\mbox{co}}\left(X \cup \{a\}\right)$ with $\mbox{dist}(x_n,A) \to 0.$
\end{itemize}
Suppose $f(a) < m$. Then there exists $\delta > 0$ such that for every $y \in B(a,\delta)$, $f(y) < m$. But (b) implies that $f(y) > m$ for all $y \in \overline{\mbox{co}}\left(X \cup \{a\}\right)$ with $y \ne a$, which yields a contradiction. Hence, $f(a) = m = \inf\left\{f(y) : y \in \overline{\mbox{co}}\left(X \cup \{a\}\right)\right\}$.\\
Take $Y = \overline{\mbox{co}}\left(X \cup \{a\}\right)$. Clearly, (a) yields that $h(X,Y) \le r$.\\
For $n \in \mathbb{N}$, let $y_n \in L_{Y,f}(1/n)$. Then $f(y_n) \le m + \frac{1}{n}$. Applying Lemma \ref{Lemma_dist} yields that $\lim_{n \to \infty}\mbox{dist}(y_n,A) = 0$.\\
Since $y_n \in L_{Y,f}(1/n)$ we also have that $\mbox{dist}(y_n,Y) \le 1/n$ and therefore, there exists $x_n \in Y$ such that $\lim_{n \to \infty}d(x_n,y_n) = 0$. Since
\[\mbox{dist}(x_n,A) \le d(x_n,y_n) + \mbox{dist}(y_n,A),\]
by (c), it follows that $x_n \to a$ and so $y_n \to a$. This implies that $\lim_{n \to \infty}\mbox{diam}L_{Y,f}(1/n) = 0$ which means that $\min(Y,f)$ is well-posed. Consequently, the proof of the denseness result is complete.

Let us now prove that $\mathcal{W}_{min}$ is a $G_\delta$ subset of $P_{b,cl,cv}(E)$. It is clear that
\[\mathcal{W}_{min} = \bigcap_{p \in \mathbb{N}}\mathcal{L}_p.\]
Thus, it suffices to prove that $\mathcal{L}_p$ is open in $P_{b,cl,cv}(E)$ for every $p \in \mathbb{N}$. Let $p \in \mathbb{N}$ and consider $X \in \mathcal{L}_p$. Pick $\sigma > 0$ such that $\mbox{diam}L_{X,f}(\sigma) < 1/p.$ Then there exists $x_0 \in X$ such that
\[f(x_0) < \inf_{x \in X}f(x) + \frac{\sigma}{3}.\]
Since $f$ is continuous, there exists $\delta < (2\sigma)/3$ such that
\[f(y) < f(x_0) + \frac{\sigma}{3} \mbox{ for every } y \in E \mbox{ with } d(x_0,y) < \delta.\]
Let $Y \in P_{b,cl,cv}(E)$ with $h(X,Y) < \delta$. We show that $Y \in \mathcal{L}_p$. Let $y \in L_{Y,f}(\sigma/3)$. Then
\[f(y) \le \inf_{z \in Y}f(z) + \frac{\sigma}{3} \mbox{ and } \mbox{dist}(y,Y) \le \sigma/3.\]
Since $h(X,Y) < \delta$, there exists $y_0 \in Y$ such that $d(x_0,y_0) < \delta$ and so 
\[f(y_0) < f(x_0) + \frac{\sigma}{3} < \inf_{x \in X}f(x) + \frac{2\sigma}{3}.\]
Thus, $f(y) < \inf_{x \in X}f(x) + \sigma$. Likewise, 
\[\mbox{dist}(y,X) \le \mbox{dist}(y,Y) + h(X,Y) < \frac{\sigma}{3} + \frac{2\sigma}{3} = \sigma.\]
Therefore, $y \in L_{X,f}(\sigma)$ and the conclusion follows.
\end{proof}

\begin{remark}
If $f$ is a continuous function, then the problem $\min(A,f)$ is well-posed in the sense of Revalski if and only if it is well-posed in the sense of Hadamard (see \cite{Revalski} for definition and proof). Hence, in the above result we can substitute the well-posedness in the sense of Revalski by the one in the sense of Hadamard.
\end{remark}

Theorem \ref{Th_f_wellp} is not only interesting by itself, but it is also important because many best approximation results follow as simple consequences thereof. For example, one can derive the following extension of a result proved in \cite{DeBlasi1}. For $X \in P_{b,cl,cv}(E)$ and $y \in E$, we denote the nearest point problem of $y$ to $X$ by $\min(X,y)$. A sequence $(x_n)_{n \in \mathbb{N}}$ in $X$  such that $d(x_n,y) \to \mbox{dist}(y,X)$ is called a minimizing sequence. The problem $\min(X,y)$ is said to be well-posed if it has a unique solution $x_y \in X$ and every minimizing sequence converges to $x_y$.

\begin{corollary}
Let $E$ be a complete geodesic space with convex metric and let $y \in E$. Then the set
\[\mathcal{W}_{min} = \left\{X \in P_{b,cl,cv}(E) : \min(X,y) \mbox{ is well-posed} \right\}\]
is a dense $G_\delta$-set in $P_{b,cl,cv}(E)$.
\end{corollary}
\begin{proof}
Take $f: E \to \mathbb{R}, f(x) = d(x,y)$ for $x \in E$ in Theorem \ref{Th_f_wellp}.
\end{proof}

\section{Acknowledgments}
We would like to thank Art Kirk for a preliminary discussion on this topic and Genaro L\' opez for bringing Georgiev's paper \cite{Georgiev} to our attention.\\
The research of the first author was partially supported by DGES, Grant MTM2009-10696-C02-01 and Junta de Andaluc\'{\i}a, Grant FQM-127. The second author was supported by programs co-financed by The Sectoral Operational Programme Human Resources Development, Contract POS DRU 6/1.5/S/3 - ``Doctoral studies: through science towards society''. She would also like to express her appreciation to the Department of Mathematical Analysis and to the Institute of Mathematics of the Univeristy of Seville (IMUS) for their support.

\end{document}